\documentclass[12pt]{paper}
 \usepackage{epsfig}
\usepackage{graphics}
 \usepackage{color}

\usepackage{color}
\usepackage{epsfig,cmmib57}
\usepackage{amssymb,amsmath,exscale}
 \numberwithin{equation}{section}

\newcommand{\ZEP}{\epsilon}

\newcommand{\ZSUno}{\sum _{n=1}^{+\ZIN}}

\newcommand{\zg}{\gamma}

\newcommand{\intp}{\int_0^{\pi}}

\newcommand{\intt}{\int_0^t}
\newcommand{\ints}{\int_0^s}
\newcommand{\intr}{\int_0^r}

\newtheorem{Theorem}{Theorem}

\newtheorem{Lemma}[Theorem]{Lemma}

\newtheorem{Remark}[Theorem]{Remark}

\newcommand{\zdiaform}{\mbox{~~\zdia}}

\newcommand{\ZOM}{\omega}
\newcommand{\zaa}{\alpha}

\newcommand{\ZDE}{\delta}
\newcommand{\zt}{\tau}
\newcommand{\zdia}{~~\rule{1mm}{2mm}\par\medskip}
\newcommand{\zthe}{\theta}

\newcommand{\ZLA}{\label}
\newcommand{\ZIN}{\infty}
\newcommand{\zProof}{{\noindent\bf\underbar{Proof}.}\ }
\newcommand{\zreal}{\Re{\textstyle e}\,}

\newcommand{\ZBI}{\bibitem}
\newcommand{\ZD}{\;\mbox{\rm d}}
\newcommand{\zl}{\lambda}
\newcommand{\ZSI}{\sigma}

\author{
S. Avdonin\thanks{Department of Mathematics
and Statistics,
   University of Alaska,
    Fairbanks, AK 99775-6660, USA (s.avdonin@alaska.edu)}
and
L. Pandolfi\thanks{Dipartimento di Scienze Matematiche ``Giuseppe Luigi Lagrange'', Politecnico di Torino, Corso Duca degli Abruzzi 24, 10129 Torino, Italy (luciano.pandolfi@polito.it)}
}
\title{A linear algorithm for the identification of a weakly singular relaxation kernel using  two boundary measurements
}

\begin{document}

\maketitle
{\bf\underline{Abstract}:}  
We consider a distributed system 
of a type which is   encountered   in the study of diffusion processes with memory and in viscoelasticity. The key feature of such system is the persistence in the future of the past actions due the memory described via a certain  \emph{relaxation kernel,} see below. 
   The parameters of the kernel have  to be inferred from experimental measurements.  Our main result in this paper is  that by using two boundary measurements the identification of a relaxation kernel which is  linear combination of \emph{Abel kernels} (as often assumed in applications) can be reduced to the solution of a \emph{(linear) deconvolution problem.}

 \medskip


\section{Introduction}

The following equation
\begin{equation}\ZLA{eq:sistema}
\zthe'(t)=\intt N(t-s)\zthe_{xx}(s)\ZD s
\end{equation}
is encountered for example in the study of diffusion processes, in particular thermodynamics with memory or nonfickian diffusion, and viscoelasticity. Here $\zthe=\zthe(x,t)$ where $t>0$ (prime denotes time derivative) and $x\in (a,b)$ (the 1-D case is enough for the problem studied in this paper).   Eq.~(\ref{eq:sistema}) represents the diffusion in a homogeneous isotropic medium in the form of a slab and $N(t)$ is usually called the \emph{relaxation kernel.}

 In the contest of diffusion processes, Eq.~(\ref{eq:sistema}) is obtained from conservation of energy when the usual Fourier law $q(x,t)=- \theta_x(x,t) $   for the flux (which leads to the standard heat or Fick equations) is modified and takes into account that the flux reacts ``slowly in time'' to the variation of 
the temperature (see~\cite{PandLIBRO,PandIDENT}),
\begin{equation}
\ZLA{eq:flusso}
q(x,t)=-\intt N(t-s)\zthe_x(x,s)\ZD s\,.
\end{equation}

It is a fact that certain physical considerations impose general restrictions to the \emph{relaxation kernel} $N(t)$. In particular, dissipativity on closed cycles implies that $N(t)$ is a   positive decreasing function and, by the second principle of thermodynamics,  $\hat N(\zl)$ (the Laplace transform) is defined (at least) for  $\zreal\zl>0$ and $\zreal \hat N(\zl)\geq 0$ when $\zreal \zl>0$ (see~\cite{Day,Gentili}). In the parliance of system theory, $\hat N(\zl)$ is a \emph{positive real}  transfer function.

Provided these general restrictions are satisfied, the kernel $N(t)$ has to be inferred from experiments taken on a sample of the material and it turns out   that (see~\cite{BykovMILLISECONDS}) when the material is homogeneous and isotropic $N(t)$ does not depend either on $x$ or the shape of the sample. Hence,  it can be identified using a slab of the material (this is the reason why we confine ourselves to one space dimension).

In several concrete cases, $N(t)$ is smooth for $t\geq 0$ and it is well approximated by a finite Prony sum but it turns out that often the memory   is largely predominant for small times and a combination of Abel kernels is more suitable to fit experimental data (see~\cite{Pipkin}):
\begin{equation}
\ZLA{eq:ABELkernels}
N(t)=\sum _{m=0}^N \frac{\zaa_m}{t^{\gamma_m}} \,,\quad \zaa_m>0\,,\ 0\leq \gamma_m<1\,.
\end{equation}

The specific function $N(t)$ (i.e. the values of the parameters $\zaa_m$ and $\zg_m$) has to be inferred by suitable measurements. Several identification algorithms have been proposed (see section~\ref{sect:Previousreferences}). Here we show that,  using two different measurements of the boundary flux, the identification problem can be reduced to a standard \emph{deconvolution problem,} i.e. to the solution of a Volterra integral equation of the first kind.

Note that 
\begin{equation}\ZLA{eq:LAPLAaBEL}
\mbox{the Laplace transform of $t^{-\zg}$ is $\Gamma(1-\gamma)\frac{1}{\zl^{1-\zg}}$}
\end{equation}
 (where $\Gamma(\zl)$ is the Euler
 Gamma function). So, the Laplace transform of the function $N(t)$ in~(\ref{eq:ABELkernels}) is a linear combinations of powers $1/\zl^\gamma$ where $\gamma\in (0,1)$. Instead, if $N(t)\in W^{{1,1}}([0,+\ZIN))$ and it is bounded, its Laplace transform is   $N(0)/\zl+o(1/\zl)$ (for $\zl\to+\ZIN$). 
These considerations suggest the following assumption, which is satisfied in particular by the kernels~(\ref{eq:ABELkernels}):
\medskip

{\bf\underline{Assumption}:}  
 \begin{itemize}
 \item[\bf a)] For every $T>0$ we have $N\in L^1(0,T)$; the Laplace transform $\hat N(\zl)$ of $N(t)$ is defined on $\Pi_0=\{\zreal\zl>0\}$ and it is a \emph{positive real transfer function } i.e. it is holomorphic on $\Pi_0$   and it transforms $\Pi_0$ into itself;
 \item[\bf b)]
  The following equality holds in $\Pi_0$:
\begin{equation}
 \ZLA{eq:expreNhat}
 \hat N(\zl)=\frac{\zaa}{\zl^\zg	} + \frac{1}{\zl^{\zg+\ZSI}} N_1(\zl)
 \end{equation} 
 where
 \begin{itemize}
 \item[\bf b1)] $\zaa>0$, $\zg\in (0,1)$;
 \item[\bf b2)]  $\ZSI>0$ and $N_1(\zl)$ is holomorphic   bounded  on $\zreal\zl>0$;
  \item[\bf b3)] we have $\zreal N_1(\zl)\zl^{- \ZSI }\geq 0$ in $\Pi_0$.
 \end{itemize}
 \item[\bf c)] There exists a polynomial $p(t)$ such that
 \[
 \intt |N(s)|\ZD s< p(t)\qquad t>0\,.
 \]
 
 \end{itemize}
 These assumptions are satisfied by the relaxation kernel~(\ref{eq:ABELkernels}).
 
  In typical applications, the imaginary part of $\hat N(\zl)$ and that of $\zl\in \Pi_+=\{\zreal\zl>0\}$ have opposite signs (see~\cite{Gentili}). This last property however is not used in the reconstruction algorithm.
 
 \begin{Remark}
 {\rm
 We note the following facts:
 \begin{itemize}
 \item
 We shall condider Eq.~(\ref{eq:sistema}) with the additional conditions
 \begin{equation}\ZLA{eq:IniboundCONDI}
 \zthe(0,t)=0\,,\ \zthe (\pi,t)=0\,,\quad \zthe(x,0)=\xi\in L^2(0,\pi)\,.
 \end{equation}
 Eq.~(\ref{eq:sistema}) with conditions~(\ref{eq:IniboundCONDI})  has been studied by several authors, see for example~\cite{Pruess}. We shall also consider the boundary condition $\zthe(0,t)= f(t) $ (and $\zthe(\pi,t)=0$). When $f$ is smooth, as  it is sufficient for the identification process described below,  the problem can be reduced  to a problem with distributed affine term, as studied in~\cite{Pruess}.

  \item
  The representation~(\ref{eq:expreNhat}) extracts as ``dominant part'' of $\hat N(\zl)$ the infinite of \emph{lower order} (instead of that of higher order). The reason is in   the transform
  (\ref{eq:LAPLAaBEL}): an infinite of higher order in time corresponds to an  infinite  of lower order in the frequency domain.

 \end{itemize}
 }
 \end{Remark}

\section{\ZLA{sect:Previousreferences}Comments and references}

Distributed systems with persistent memory are widely used in applied sciences. As stated already, applications range from the analysis of diffusion processes in the presence of complex molecular structure~(see~\cite{DEkEELIUHinestroza}) to thermodynamics and viscoelasticity, see~\cite{cristensen,FabrizioOWENS,Kolsky}.
 They are characterized by the fact that past action continues to affect the system, and the past effect is described via a convolution integral. So, the kernel of the convolution integral is a material property of the process, and different processes are best described with kernels in suitable classes. When the system has a ``short'' memory,  the behavior of the system is best represented with a kernel whose values in a short time interval $(0,\ZEP)$ dominate, as described in~\cite{Pipkin}, and quite often a kernel which is the sum of a smooth function and of an Abel kind kernel is used, see~\cite{Dinzart,ATANA}. 
 The memory kernel being a material property of the process, it has to be identified as the result of suitable experimental measurments and much work has been devoted to this problem both by engineers and by mathematicians. The idea is to excite the body with suitable signals and to observe a suitable output. In the engineering literature it is often assumed that the kernel depends on a small number of parameters and the body is excited with a ``simple'' signal, for example a constant deformation suddenly applied at time $t=0$. The resulting stress is then observed and compared with the stress theoretically computed, so to find the values of the parameters which ``best'' fit the measures (according a a certain index, usually a quadratic index, 
 see~\cite{Gerlac,GolubKozbarRegulina}). More refined methods use periodic signals and in essence compute the frequency response of the system (in a certain frequency range) from which the memory kernel is  recovered  via inverse Laplace transform, see~\cite{cristensen,Kolsky}. Mathematical papers use a different turn of ideas: 
a suitable output $y(t)$ is associated to the system and 
  $y'(t)$ is computed.  This gives    an additional  differential equation for $N(t)$ (which depends on the initial and boundary conditions, which should be known and have suitable properties).  In this way we get a system of nonlinear Volterra integrodifferential equations in the unknown $(\zthe,N)$, which is at least locally solvable (for small $t$), and globally solvable under suitable assumptions.   The solution of this system of equations   gives   both $N(t)$, which is the object of interest, and the solution $\zthe(t)$, which depends on the chosen initial and boundary conditions. See~\cite{FabrizioOWENS,GRASSEKABANlore,GRASSEKABAN,Janno,Lerenzilibro}.

 This method is highly non linear. Assume that it has been applied to an interval $[0,\tau]$ so that the pair $(N(t),\zthe(t))$ is known for every $t\in [0,\tau]$. As noted 
  in~\cite{Guidetti} , due to the convolutional structure of the equation, the identification of $(N(t),\zthe(t) )$ on $(\zt,T]$ can then be reduced to a linear deconvolution problem.

  The method we present here reduces the identification of the relaxation kernel to a linear deconvolution problem and the parallel identification of a specific evolution $\zthe(t)$ is not needed. Note that this is an important simplification because the object of interest is a real valued function $N(t)$ and $\zthe(t)=\zthe(\cdot,t)$ is an $L^2$-valued function of time, hence far   more  complex then the required quantity $N(t)$. This method has already been justified in~\cite{PandIDENT,PandMULTI} when $N(t)$ is smooth for $t\geq 0$ but the time domain proofs of these papers seems not extendable to the case of singular kernels, and   so in this paper     the algorithm is justified by using frequency domain techniques.

\section{The algorithm}

The algorithm we propose requires the computation of the flux through  the boundary, say at the end $x=\pi$, due to suitable excitations of the sample:

\begin{enumerate}
\item The first measurement  keeps both the end of the sample at fixed temperature (equal zero without restriction) and we measure  the flux     (denoted $q_{\xi}(\pi,t)$) through  the point  $x=\pi$ due to a suitable initial temperature $\xi$.   So we solve  Eq.~(\ref{eq:sistema}) with the following conditions
\begin{equation}\ZLA{eq:condiINI}
\zthe(x,0)=\xi(x)\,,\qquad \zthe(0,t)=\zthe(\pi,t)=0\,.
\end{equation} 
\item
 In a second measurement  we keep the initial temperature and the temperature at $x=\pi$ at a fixed value (equal zero without restriction), while the temperature at the left end $x=0$ is regulated, i.e. we solve Eq.~(\ref{eq:sistema}) with the following conditions
\begin{equation}
\ZLA{eq:CondiPRIMAmisura}
\zthe(x,0)=0\,,\qquad \zthe(0,t)=f(t)\,,\quad \zthe(\pi,t)=0\,.
\end{equation} 
We measure the resulting flux    through   the point  $x=\pi$.
We denote   $q^f(\pi,t)$ this flux.
\end{enumerate}
We shall see that $N(t)$, $t\in [0,T]$, can be identified from these measurements taken on the interval $  (0,T)$. What is more important is the fact that \emph{the identification is reduced to the solvability of a standard  deconvolution problem.}

Some comments are in order:
\begin{enumerate}
\item the algorithm extends to the case of the singular kernel $N(t)$ in~(\ref{eq:ABELkernels}) the algorithm 
introduced in~\cite{PandIDENT,PandMULTI}  when $N(t)\in C^3([0,T])$. 
 The proofs given in these papers heavily depend on smoothness of the relaxation kernel even for $t=0$. So, the contribution of the present paper is the justification of the algorithm in the case of weakly singular kernels.
\item The time varying boundary temperature $f(t)$ can be easily imposed, provided that it is smooth and $f(0)=0$ (consistent with the condition $\zthe(x,0)\equiv 0$, see~(\ref{eq:CondiPRIMAmisura})) and in practice, bounded. So, $f(t)$ has the representation
\[f(t)=\intt g(s)\ZD s \]
 and without restriction we can assume $g \in C [0,T]$.
\item In practice, it is difficult to impose an initial temperature $\xi(x)$, unless $\xi(x)$ is quite special. We shall see that   the special temperature needed in the identification algorithm can be easily imposed.
\end{enumerate}

Finally, we comment again on the expression of the flux. The (density of the) flux at time $t$ and position $x$ 
   is given by~(\ref{eq:flusso}) so that 
\[
  q(\pi,t)=-\intt N(t-s)\zthe_x(\pi,s)\ZD s\,.
\]
It is not at all obvious that the expression of $q(\pi,t)$ makes sense, and in what kind of function or distribution space. This point is clarified in the Appendix.

  \section{\ZLA{sect:fluxSOLUdefi}The solution and the flux}

The fact that $\hat N(\zl)$ is a positive real transfer function  implies existence of the solution when the boundary input $f$ is equal zero (the case $f\neq 0$ will be examined below). To see this, let us introduce in $L^2(0,\pi)$ the operator $A=\Delta$ with domain $\{\phi\in H^2(0,\pi)\,,\ \phi(0)=\phi(\pi)=0\}$ and rewrite Eq.~(\ref{eq:sistema}) in the form
\begin{equation}\ZLA{eq:sisteMODIF}
\zthe(t)=A\intt a(t-s)\zthe(s)\ZD s+ F(x,t)\,,\quad  a(t)=\intt N(s)\ZD s\,. 
\end{equation}
Here $F(x,t)= \zthe(x,0)$ (constant in time).

For every $F\in C([0,T];L^2(0,\pi))$
there exists a unique function $\zthe(t)\in C([0,T];L^2(0,\pi))$ such that $A\intt a(t-s)\zthe(s)\ZD s\in C([0,T];L^2(0,\pi))$ and which solves~(\ref{eq:sisteMODIF}) (for every $T>0$). The solution depends continuously on the affine term $F$. This follows from~\cite[Propositions~1.2, 1.3 and Corollary~1.2(i)  Chapt.~1]{Pruess}  since $\hat N(\zl)$ is a positive real transfer function.

In the case the boundary input $f\neq 0$ we shall prove:
\begin{Theorem}
Let $f\in C^1([0,T])$ and $f(0)=0 $.
Then Eq.~(\ref{eq:sistema}) with  the initial and 
 boundary conditions~(\ref{eq:CondiPRIMAmisura}) admits a unique solution $\zthe\in C([0,T];L^2(0,\pi))$.
 \end{Theorem}
Of course, when $\zthe(x,0)=\xi(x)\neq 0$ the previous results can be combined.
\zProof We proceed formally and then we justify the results.  For every $t$ we introduce $u_0(x,t)= \frac{\pi-x}{\pi}f(t)$ and we note that $\Delta u_0(x,t)=0$ for every $t$. Let
\[
\xi(x,t)=\zthe(x,t)-u_0(x,t)
\]
(so that $\xi(x,0)=0$ since $f(0)=0$).
If $\zthe$ should be a solution of Eq.~(\ref{eq:sistema}) with  conditions~(\ref{eq:CondiPRIMAmisura}), then $\xi$ should solve
\[
\xi'=\intt N(t-s) A\xi(s)\ZD s - \frac{\pi-x}{\pi}f'(t) 
\]
 with zero boundary conditions. 
(As above, $A$ is the laplacian, restricted to functions which are zero at $x=0$ and $x=\pi$).
 The function $f'(t)$ is continuous    so that, as already noted, the existence of the solution 
$\xi(x,t)$ (in the mild sense specified above)   again  follows  from~\cite[Propositions~1.2, 1.3 and Corollary~1.2(i)  Chapt.~1]{Pruess}.

The solution $\zthe=\zthe(x,t)$ is by definition the function $\zthe(x,t)=\xi(x,t)+((\pi-x)/\pi) f(t)$.
\zdia

\begin{Remark}
The condition that $f$ is   continuously   differentiable is sufficient  for the identification problem. So, we don't invest time to weaken this assumption.\zdia
\end{Remark}

Now we proceed also by separation of variables and we find formulas which are used in order to describe and justify the reconstruction algorithm. Let
\[
\phi_n=\sqrt{\frac{2}{\pi}}\sin nx\,,\quad n\in \mathbb{N}, \quad \mbox{so that}\quad  \frac{\ZD^2}{\ZD x^2} \phi_n(x)=-n^2 \phi_n(x)\,.
\]
It is known from the theory of Fourier series that $\{\phi_n\}$ is an orthonormal basis of $L^2(0,\pi)$.
We solve Eq.~(\ref{eq:sistema})   with   the initial condition $\zthe(x,0)=\xi(x)$ and the boundary conditions $\zthe(0,t)=f(t)$, $\zthe(\pi,t)=0$.
We expand the initial condition $\xi(x)$ and the solution $\zthe(x,t)$ of Eq.~(\ref{eq:sistema})  in   the  sine series:
\[
 \xi(x)=\frac{2}{\pi}\ZSUno \xi_n\sin nx\,,\qquad \zthe(x,t)=\frac{2}{\pi}\ZSUno \zthe_n(t)\sin nx 
\]
where, for every $n\in \mathbb{N}$,  
\[
\xi_n:=\intp \xi(x)\sin nx\ZD x\,,\qquad \zthe_n(t):=\intp \zthe(x,t)\sin nx\ZD x \,.
\]
Note that $\zthe$ can be expanded in the sine series   since $\zthe\in C([0,T];L^2(0,\pi))$.

The functions  $\zthe_n(t)$ solve 
\[
\zthe'_n=-n^2\intt N(t-s) \zthe_n(s)\ZD s+n\intt N(t-s) f(s)\ZD s\,,\qquad \zthe_n(0)=\xi_n\,.
\]
This is a scalar integrodifferential equation whose solution can be represented using the \emph{variation of constants  formula} (similar to~\cite[Formula~(5.7)]{PandLIBRO}. In this reference $N(t)$ is regular but the derivation of the variation of constants formula does not use this condition).
We introduce the solution $z_n(t)$ of
\begin{equation}\ZLA{equaDIz}
z_n'(t)=-n^2\intt N(t-s) z_n(s)\ZD s\,,\qquad z_n(0)=1\,.
\end{equation}
Then we have
\begin{align} 
\nonumber \zthe_n(t)&=  \xi_n z_n(t) +  \left [ n \intt z_n(t-s)\ints N(s-r) f(r)\ZD r\,\ZD s\right ]=\\
\ZLA{eq:ExpaCompoSeRIE} &=\xi_n z_n(t) -\frac{1}{n} \intt f(r)z_n'(t-r)\ZD r 
\end{align}
(note that the prime denotes derivatives   with respect to the variable $t$) and

\begin{align}
\nonumber \zthe(x,t)&=\frac{2}{\pi}\ZSUno\left (\sin nx\right ) \xi_n z_n(t)-\\
\ZLA{eq:FORMAditheta}
& -\frac{2}{\pi}\ZSUno 
\left (\sin nx\right )\left [
\frac{1}{n} \intt f(r) z_n'(t-r)\ZD r
\right ]
\end{align}
 
We conclude this section with the following  observation (see also~\cite[~ Corollary~1.2  Chapt.~1]{Pruess}), which will be used in the Appendix.
The assumption that $\hat N(\zl)$ is a  positive real transfer function implies $|z_n(t)|\leq 1$ for every $n$ and every $t>0$.   Indeed,  multiplying   both the sides of~(\ref{equaDIz}) by $z_n(t)$ and integrating we get
\[
z_n^2(t)=1-n^2 \intt z_n(s)\int_0^s N(s-r) z_n(r)\ZD r\,\ZD s\leq 1 
\]
  and note that   the iterated convolution is non    nonpositive  because $\hat N(\zl)$ is a positive real transfer function (see~\cite{Pruess}).

\subsection{The flux}

Now we compute the flux through   the point  $x=\pi$.
We proceed formally and we find certain  series whose convergence will be justified in the appendix.

  At first,  we consider  the case $f=0$ and $\xi\neq 0$. The flux in this case is denoted $q_\xi(x,t)$.   We have, by using~(\ref{eq:flusso})  and the series~(\ref{eq:FORMAditheta})

\begin{align}
\nonumber\frac{\pi}{2}q_\xi(x,t)=-\intt N(t-s)\left (
\ZSUno \left (n\cos nx\right ) \xi_n z_n(t)
\right )\ZD s =\\
-\ZLA{eq:FluxDATINIZ}\ZSUno \cos nx\left (
n\intt N(t-s) z_n(s)\ZD s
\right )\xi_n  = \ZSUno \cos nx \frac{1}{n} z'_n(t) \xi_n \,.
\end{align}
Let us consider now the flux $q^f(x,t)$ due to the boundary temperature $f$ when $\xi=0$.  We assumed  that 
\[
f(t)= \intt g(s)\ZD s, \ g \in C[0,T].
\] 
 Then we have  
\begin{align*}\nonumber
&\frac{\pi}{2}q^f(x,t)=\intt N(\zt)\left (\ZSUno \left (\cos nx\right ) \int_0^{t-\zt} f(r)z_n'(t-\zt-r)\ZD r\right )=
\\ \nonumber
&-\intt N(\zt)\left (
\ZSUno \left (\cos nx\right )\int_0^{t-\zt}\left (\frac{\ZD}{\ZD r}z_n (t-\zt-r)\right )\intr g(\nu)\ZD\nu \ZD r 
\right )\ZD\zt=
 \\ \nonumber
& -\left (\intt N(\zt)\int_0^{t-\zt}g(\nu)\ZD\nu \ZD \zt\right ) \ZSUno \cos nx+\\
\nonumber &+
\nonumber \intt N(\zt)\left (\int_0^{t-\zt} g(r)\ZSUno \left (\cos nx\right ) z_n(t-\zt-r)\ZD r \right )\ZD\zt\,.
\end{align*}
We recall   that 
 \[
\ZSUno\cos nx=\pi\ZDE(x)-\frac{1}{2}
\]
($\ZDE(x)$ is the Dirac delta and the convergence is in the sense of the distributions). 
  Therefore,  
\begin{align*}
&\frac{\pi}{2}q^f(x,t )=\left (\pi\ZDE(x)\right )\intt N (\zt)f(t-\zt)\ZD \zt-\\
&-\frac{1}{2} \intt N (\zt)f(t-\zt)\ZD \zt-
\intt g(r)\ZSUno \frac{1}{n^2}z'_n(t-r)\cos nx \ZD r\,.
\end{align*}

We are interested in the flux through the right hand $x=\pi$. 
So we measure
\begin{align}
\ZLA{flussoxi}&y_\xi(t)=\frac{\pi}{2}q_\xi(\pi,t)=\ZSUno (-1)^n\frac{1}{n}\xi_n z_n'(t)\,,\\
\ZLA{flussof}& Y^f(t)=\frac{\pi}{2}q^f(\pi,t)=-\frac{1}{2}\intt \tilde N(t-\zt)g(\zt)\ZD\zt-\intt g(r) y_{\xi_0}(t-r)\ZD r
\end{align}
  where
\begin{equation}\ZLA{formaDATINIZ}
\tilde N(t)=\intt N(s)\ZD s\,,\qquad 
\xi_0(x)=\frac{2}{\pi} \ZSUno \frac{1}{n} \sin nx=\frac{1}{2}(\pi-x)\,.
\end{equation}

Note that these are formal expressions   for $y_\xi$ and $Y^f$  to be justified below in the appendix. Granted these formulas, we observe:
\begin{itemize}
\item the initial temperature $\xi_0=(\pi-x)/2$ is in the form of a ramp. When the system is ``stable'', i.e. the free evolution tends to a stationary temperature  as usually happens in applications, an initial condition in the form of a ramp    is easily achievable as the equilibrium temperature when the ends of the bar are kept at  constant temperatures  for a time long enough, see~\cite{PandMULTI}. 
 
\item  once the measurement $y_{\xi_0}(t)$ has been obtained, the determination of $\tilde N (t) $ from~(\ref{flussof}) is a standard deconvolution problem, which can easily be solved by existing methods. If it is known, or postulated, that the relaxation kernel has the form~(\ref{eq:ABELkernels}) then the  parameters $\zaa_m$ and $\zg_m$ can be computed from $\tilde N(t)$. The determination of $N(t)$, i.e. the computation of the numerical derivative of $\tilde N(t)$, is not needed.

\end{itemize}

This is the identification algorithm, which extends to weakly singular kernels the algorithm introduced   in~\cite{PandIDENT,PandMULTI} for smooth  kernels.   Of course, the algorithm is not yet justified since we must still prove convergence of the series~(\ref{flussoxi}) and~(\ref{flussof}). The proof is in the appendix.

\section{Appendix: justification of the formulas}

   Here we justify the formulas for the flux, and so the correctness of the algorithm we proposed. 
In fact, we prove the following result:
\begin{Theorem}
Let $T>0$. We have:
\begin{enumerate}
\item $y_\xi\in L^2(0,T)$ for every $\xi\in L^2(0,\pi)$ and the transformation $\xi\mapsto y_\xi$ is linear and continuous from $L^2(0,\pi)$ to $L^2(0,T)$;
\item $Y^f\in L^2(0,T)$  for every $f\in H^1_-(0,T)=\{ f \in H^1(0,\pi)\,,\ f(0)=0\}$ and the transformation $f\mapsto Y^f$ is continuous from $H^1_-(0,T)$ to $L^2(0,T)$. 
\end{enumerate}
\end{Theorem}

The second statement follows easily from the first one, as   seen from formula~(\ref{flussof}). Using the first statement, i.e.  $y_{\xi_0}\in L^2(0,T) $, and Young inequality it is easily seen that $Y^f\in L^2(0,T)$ and that it depends continuously on   $g\in L^2(0,T)$, hence   on   $f\in H^1_-(0,T)$.

So, it is enough that we prove the first statement. The proof is in the frequency domain, via   the   Laplace transform and Parseval identity.

We must prove the convergence in $L^2(0,T)$ of the series~(\ref{flussoxi}) where $\{\xi_n\}\in l^2$. So it is sufficient to prove boundedness in $L^2(0,T)$ of the sequence
\[
\left \{
\frac{z'_n(t)}{n}
\right \}=\left \{
-n\int_0^t N(t-s)z_n(s)ds
\right \}\,.
\]

   We know already that   $ |z_n(t)|<1 $ so that, for every fixed $n$,      Assumption {\bf c)} shows that
  \[ 
  | z_n'(t) |\leq n^2\intt |N(s)|\ZD s\quad \mbox{is polynomially bounded.}   
   \]
    Hence, for every $\beta>0$, $e^{-\beta t}z_n'(t)\in L^2(0,+\ZIN)$. We prove that for a suitable value of $\beta$ there exists $M>0$ \emph{independent of $n$} and such that  
   \[
 \int_0^{T} e^{-2\beta t} \left ( \frac{z_n'(t)}{n} \right )^2\ZD t \leq
 \int_0^{+\ZIN} e^{-2\beta t} \left ( \frac{z_n'(t)}{n} \right )^2\ZD t <M   
   \]
   which implies boundedness in the standard $L^2(0,T)$-norm:
   \[
    \int_0^{T}   \left ( \frac{z_n'(t)}{n} \right )^2\ZD t \leq M
   \]
   (for a different constant  $M$).
   
  Let us consider now any function $h$ defined on $(0,+\ZIN)$ and let $\hat{}$ denote the Laplace transform. Note that if $\hat h(\zl)$ exists, then the Laplace transform of $e^{-\beta t}h(t)$ is $\hat h(\beta+\zl)$ and  the following result hold (see~\cite[Ch.~6]{Koosis}):
   \begin{Theorem}
Let $\beta\geq 0$.  The function $e^{-\beta t}h(t)$ belongs to $L^2(0,+\ZIN)$ if and only if $\hat h(\zl)$ is holomorphic in
   \[
  \Pi_\beta= \{\zl\,: \ \zreal\zl>\beta\}\,,\quad  
  \mbox{and}\quad 
\sup _{\xi>0} \int_{-\ZIN}^{+\ZIN} |\hat h(\xi+i\ZOM)|^2\ZD \ZOM<+\ZIN\,.
\]
Furthermore, under these conditions
\[
\pi \int_0^{+\ZIN} e^{-2\beta t} h^2(t)\ZD t=\int _{-\ZIN}^{+\ZIN} |\hat h(\beta+i\ZOM)|^2\ZD \ZOM\,.
\]
   \end{Theorem}
   
   We apply this result to the functions
   \[
h(t)=   \frac{z'_n(t)}{n}
 = 
-n\int_0^t N(t-s)z_n(s)ds\,,  \ \  n \in \mathbb{N}. 
   \]

 The Laplace transform of
 \[
-\frac{1}{n}z_n'(t)=  n \intt N(t-s) z_n(s)\ZD s
 \] 
is  
 
\begin{equation}
\ZLA{eq:LaplaTRANSfoFlu1}
n\hat N(\zl)\frac{1}{\zl+n^2\hat N(\zl)} =\frac{n}{n^2+\zl/\hat N(\zl)}
\,.\end{equation}
We   evaluate

 \begin{equation}
 \ZLA{inteDUEdaCalco}
 \int _{-\ZIN}^{+\ZIN} \frac{1}{\left |n^2+\zeta(\beta+i\ZOM)\right |^2}\ZD\ZOM 
 \end{equation}
where $\zeta(\beta+i\ZOM)$ is the function

 \begin{equation}\ZLA{DefiZita}
\ZOM \mapsto\zeta(\beta+i\ZOM)= \frac{\beta+i\ZOM}{\hat N(\beta+i\ZOM)}=   \frac{(\beta+i\ZOM)^{1+\gamma}}{\alpha+N_1(\beta+i\ZOM)(\beta+i\ZOM)^{-\ZSI}}\,.  
 \end{equation}
 Integral~(\ref{inteDUEdaCalco})  is an integral on the vertical line $\xi=\beta+i\ZOM$ of the complex plane.
 
   To prove (5.1) we will prove that, for a suitable $\beta>0$,  the integral  in~(\ref{inteDUEdaCalco}) decays of the order $1/n^2$.  We  consider the integral on $(0,+\ZIN)$; the integral on $(-\ZIN,0)$  can be treated analogously. 
 
 It is easier to follow the arguments by looking at the  figure. 

 \begin{center}
\includegraphics[width=10cm]{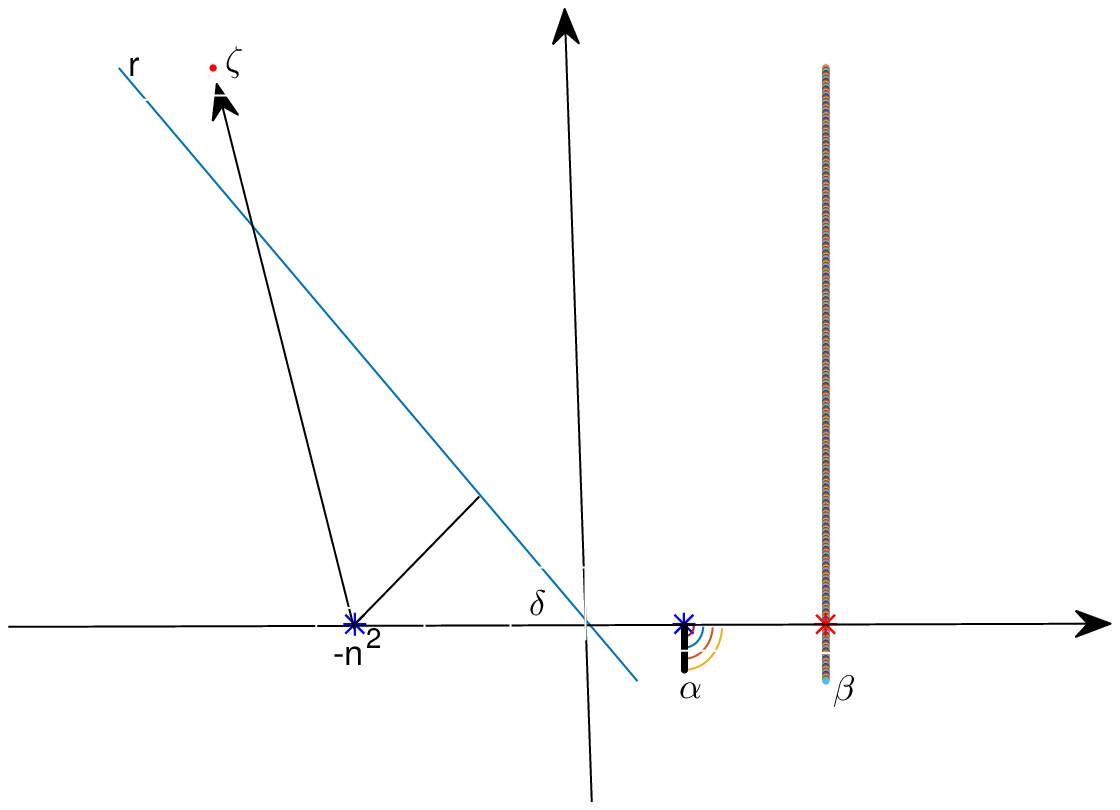} 
   \end{center}
 
 By assumption, in   the complex plane the points $\alpha+N_1(\beta+i\ZOM)(\beta+i\ZOM)^{-\ZSI}$ belong 
 to a semidisk  contained  in $ \{\zreal\zl>\zaa\}$ (and typically in   $   \mathcal{I}m\zl<0$ when $\ZOM>0$, see~\cite{Gentili}, but this property is not used).
  Boundedness of $N_1(\zl)$   implies  that the radius of this disk (whose center is $(\zaa+i0)$) is less then $M(\beta^2+\ZOM^2)^{-\ZSI/2}$ for a suitable constant $M$.
  
  Now we explain how $\beta$ has to be chosen.  The argument of $(\beta+i\ZOM)^{1+\zg}$
  is less then $(\pi/2)(1+\zg)<\pi$ since $\zg<1$.  The argument of $\alpha+N_1(\beta+i\ZOM)(\beta+i\ZOM)^{-\ZSI}$ tends to zero for $\beta\to+\ZIN$   and so there exists a value $\beta_0 $ such that if $\beta>\beta_0  $ then
   ${\rm Arg}\,\zeta(\beta+i\ZOM)<\pi-\ZDE$   for some $\delta >0$. 
   Furthermore, by increasing $\beta$, we can reduce the norm of $N_1(\beta+i\ZOM)/(\beta+i\ZOM)^\ZSI$ as much as we wish.  The value of $\beta$ is chosen so to have
\begin{equation}\ZLA{eq:condiSUbeta}
\left\{
\begin{array}{l}
{\rm Arg}\,\zeta(\beta+i\ZOM)<\pi-\ZDE\in (\frac{\pi}{2},\pi)\,,\\
\left |\frac{N_1(\beta+i\ZOM)}{(\beta+i\ZOM)^\ZSI}\right |<\zaa\,.
\end{array}
\right.
\end{equation}  
 
  This is the value of $\beta$ which is fixed in the following computation.

In the figure, the points   $\zeta(\beta+i\ZOM )$  are points of the complex plane in $\zreal\zl<0$, above the line marked $r$. The angle of $r$ with the \emph{negative} real axis is
 \[
\ZDE\in (0,\pi/2)\,.
 \]
  
 Then we have
 \[
\zeta(\beta+i\ZOM)=R(\ZOM)\left [ -\cos\left (\ZDE+\ZEP(\ZOM) \right )+i\sin\left (\ZDE+\ZEP(\ZOM) \right )\right ]\,,\quad R(\ZOM)=|\zeta(\beta+i\ZOM)|
 \] 
 and $\ZEP(\ZOM)>0$,  $\ZEP(\ZOM)\to 0$ for $\ZOM\to \ZIN$.

 From  geometrical point of view, the denominator $|  n^2+\zeta(\beta+i\ZOM)|^2$ in~(\ref{inteDUEdaCalco}) 
 is the squared distance of  $\zeta(\beta+i\ZOM)$ to the point $(-  n^2+i0)$ (the length  of the arrow in the figure).
 
 The distance of $(-  n^2+i0)$   from  the line $r$ is 
 $  n^2 \sin\ZDE$, reached by the point of $r$ whose distance from the origin is $  n^2 \cos\ZDE$.
 
 This observation suggests that we decompose the integral in~(\ref{inteDUEdaCalco}) as
 \begin{equation}\ZLA{eq:decoINTE}
\int_0^{+\ZIN}=\int_0^{   n^2 \cos\ZDE}  +\int_{   n^2 \cos\ZDE}^{+\ZIN}\,.
 \end{equation}
 We estimate the two integrals separately. 
We need an asymptotic estimate which holds for large $n$. For a reason we shall see, we consider those $n$ for which
the following inequalities hold:
\begin{equation}
 \ZLA{eq:CondiSOTTOdiN}
\begin{array}{llll}
{\bf a)}& n^2>\frac{(2\zaa)^{1/\zg}}{\cos\ZDE}\,, &
{\bf b)}& 
n^2>\frac{(2\zaa)^{1+\zg}}{\cos\ZDE} \,.
\end{array}
 \end{equation}

 The following estimate for the first integral in the right hand side of~(\ref{eq:decoINTE}) does not use~(\ref{eq:CondiSOTTOdiN}):
 
 \[
\int_0^{ n^2 \cos\ZDE}  \frac{1}{|n^2+\zeta(1+i\ZOM)|^2  }\ZD \ZOM\leq  \int_0^{   n^2 \cos\ZDE} 
 \frac{1}{n^4\sin^2\ZDE}
\ZD \ZOM=  \frac{1}{n^2}\frac{\cos\ZDE}{\sin^2\ZDE}\,.
 \]
 Once multiplied with $n^2$, this term remains bounded, as we wished.
 
 We must prove an analogous property for the second integral. We note:
 
\begin{align*}
&\left |
 n^2+\zeta(1+\ZOM)
\right |^2=
 \left |
\left (n^2-R(\ZOM)\cos (\ZDE+\ZEP(\ZOM)\right )+iR(\ZOM)\sin\left (\ZDE+\ZEP(\ZOM)\right )
 \right |^2=\\
 &= n^4+R^2(\ZOM)-2  n^2R(\ZOM)\cos\left (\ZDE+\ZEP(\ZOM)\right )\geq\\
 & \geq n^4+R^2(\ZOM)-2  n^2R(\ZOM)\cos \ZDE  =
 \left (R(\ZOM)- n^2\cos\ZDE\right )^2+ n^4\sin^2\ZDE\,.
\end{align*}

We recall that $R(\ZOM)=|\zeta(\beta+i\ZOM)|$,   $\cos\ZDE>0$ and use  the condition  {\bf a) } in~(\ref{eq:CondiSOTTOdiN}).  We prove:
 
\begin{Lemma}\ZLA{eq:lemmaDISU}
The conditions $\ZOM \geq  n^2\cos\ZDE$ and~(\ref{eq:CondiSOTTOdiN})  imply the following inequalities:
 \[
\left\{\begin{array}{ll}
{\bf i)} &\left | N_1(\beta+i\ZOM)/(\beta+i\ZOM)^{\ZSI}\right |<\zaa\,,\\
{\bf ii)} & R(\ZOM)> n^2\cos\ZDE\,,\\
{\bf iii)} &  \left (R(\ZOM)-n^2\cos\ZDE\right )^2 \geq \left (\frac{1}{2\zaa}\ZOM^{1+\zg}-n^2\cos\ZDE\right )^2\,.
\end{array}\right.
 \]
 
\end{Lemma}

\zProof  
 Property  {\bf i)} is the second inequality  in~(\ref{eq:condiSUbeta}) and it is our choice of the fixed value of $\beta$.
 
 We use {\bf i)} and condition~{\bf a)} in~(\ref{eq:CondiSOTTOdiN}) to see that
 \begin{align*}
R(\ZOM)=\frac{(\beta^2+\ZOM^2)^{(1+\zg)/2}}{\left |\zaa+N_1(\beta+i\ZOM)/(\beta+i\ZOM)^\ZSI\right |}\geq \frac{(\beta^2+\ZOM^2)^{(1+\zg)/2}}{\zaa-|N_1(\beta+i\ZOM)|/(\beta+\ZOM^2)^{\ZSI/2}}\\
 \geq  \frac{(\beta^2+\ZOM^2)^{(1+\zg)/2}}{\zaa}>n^2\cos\ZDE \,.
 \end{align*}
    In fact,  by using  $\ZOM>n^2\cos\ZDE$ and inequality~{\bf a)} in~(\ref{eq:CondiSOTTOdiN}) we see that  
\begin{align*}
\frac{(\beta^2 +\ZOM^2)^{(1+\zg)/2}}{\zaa}
 \geq \frac{\left (\beta^2 +n^4\cos^2\ZDE\right )^{(1+\zg)/2}}{\zaa}
 \geq \left ( n^2\cos\ZDE\right)\frac{\left (n^2\cos\ZDE\right )^\zg}{\zaa}\geq n^2\cos\ZDE\,.
\end{align*}
I.e.   we have checked the  property {\bf ii)}.

Finally we prove~{\bf iii)}. When $\ZOM>  n^2\cos\ZDE$, property   {\bf a)}   in~(\ref{eq:CondiSOTTOdiN}) gives
\[
\frac{1}{2\zaa}\ZOM^{1+\zg}>\frac{1}{2\zaa}\left (  n^2\cos\ZDE\right )^{1+\zg}>n^2\cos\ZDE, \quad \mbox{so that } \  \frac{1}{2\zaa}\ZOM^{1+\zg}-n^2\cos\ZDE>0 
\]
and    we can confine ourselves to prove $R(\ZOM)>\frac{1}{2\zaa}\ZOM^{1+\zg}$.
This holds true because   inequality {\bf i)} shows that
$ |\zaa+N_1(\beta+i\ZOM)/(\beta+i\ZOM)^{\ZSI}|<2\zaa$ and so
\begin{align*}
R(\ZOM)=\frac{(\beta+\ZOM^2)^{(1+\zg)/2}}{|\zaa+N_1(\beta+i\ZOM)/(\beta+i\ZOM)^{\ZSI}|}>\frac{1}{2\zaa}\ZOM^{1+\zg}\,,\zdiaform
\end{align*}

In conclusion, we have
 \begin{equation}
 \ZLA{eq:stimaPERlaCODA}
\left |n^2+\zeta  (\beta+i\ZOM)\right |^2 \geq \left(\frac{1}{2\zaa}\ZOM^{1+\zg} -  n^2\cos\ZDE\right )^2+n^4\sin^2\ZDE\,.
 \end{equation}
 
  We use~(\ref{eq:stimaPERlaCODA}) in order to give an upper estimate of the 
  second integral in   the right hand side  of~(\ref{eq:decoINTE}). The upper estimate is an integral which can be computed explicitly. 
  
  In the following computation, line~(\ref{eq:STIMAfinaA}) is obtained by using the substitution $\ZOM^{1+\zg}=2\zaa s$ in the   integral at the previous line
and the inequalities between the lines~(\ref{eq:STIMAfinaA}) and~(\ref{eq:STIMAfinaB}) holds   for   $n$ so large that condition {\bf b)} in~(\ref{eq:CondiSOTTOdiN})) holds, because     in this case $s \geq (1/2\zaa)\left (n^2\cos\ZDE\right )^{1+\zg}>1$.

 \begin{align}
\nonumber & \int _{( n^2 \cos\ZDE) } ^{+\ZIN}\leq \int _{( n^2 \cos\ZDE) } ^{+\ZIN} \frac{1}{\left  (  (1/2\zaa)\ZOM^{1+\zg} - n^2\cos\ZDE\right )^2 + n^4\sin^2\ZDE}\ZD\ZOM=\\[2mm]
\ZLA{eq:STIMAfinaA} &=\frac{(2\zaa)^{1/(1+\zg)}}{1+\zg} \int _ {(1/2\zaa)( n^2 \cos\ZDE)^{1+\zg}} ^{+\ZIN}
 \frac{1}{s^{\zg/(1+\zg)}}
  \frac{1}{\left  ( s- n^2\cos\ZDE  \right )^2+ n^4\sin^2\ZDE }\ZD s\leq \\[2mm]
  \ZLA{eq:STIMAfinaB} &=\frac{(2\zaa)^{1/(1+\zg)}}{1+\zg} \int _ {(1/2\zaa)( n^2 \cos\ZDE)^{1+\zg}} ^{+\ZIN}
  \frac{1}{\left  ( s- n^2\cos\ZDE  \right )^2+ n^4\sin^2\ZDE }\ZD s\leq \\[2mm]
\nonumber &\leq \frac{(2\zaa)^{1/(1+\zg)}}{1+\zg} \int _ {-\ZIN} ^{+\ZIN}  \frac{1}{\left  ( s- n^2\cos\ZDE  \right )^2+ n^4\sin^2\ZDE }\ZD s\leq\\[2mm]
\nonumber &\leq \frac{\pi(2\zaa)^{1/(1+\zg)}}{1+\zg}\frac{1}{n^2\sin\ZDE} 
 \end{align}
 as wanted.

We sum up: we proved that the sequence $\{z_n'(t)/n\}$ is bounded in $L^2(0,T)$ for every $T$.   \emph{This implies that $y_\xi(t)\in L^2(0,T)$ for every $T>0$} (and depends continuously on $\xi\in L^2(0,\pi)$). In particular we have also $y_{\xi_0}\in L^2(0,T)$ and  \emph{this in particular  shows that $Y^f(t)$ belongs to $L^2(0,T)$} (and depends continuously on $g=f'\in L^2(0,T)$\/).  Hence the proposed algorithm is now justified.

 \section{  Acknowledgments}
 
This
work has been done in the context of a visit of  the first author to the Dipartimento di Scienze Matematiche  ``G.L. Lagrange'' of the  Politecnico
di Torino in June 2016, supported by GNAMPA-INDAM. It fits into
the research programs of GNAMPA-INDAM and of the   ``Groupement de Recherche en Contr\^ole des EDP entre la France et l'Italie (CONEDP-CNRS)''.

The research of Sergei Avdonin was supported in part by the National Science Foundation,
 grant DMS 1411564 and by the Ministry of Education and Science of Republic of Kazakhstan under the grant No. 4290/GF4.


\begin{thebibliography}{99}
\ZBI{ATANA}Teodor M. Atanackovi\'c,
Stevan Pilipovi\'c,
Bogoljub Stankovi\'c,
Dušan Zorica,
\emph{Fractional calculus with applications in Mechanics,}  John Wiley and Sons, Hoboken NJ, 2014.
\ZBI{cristensen} R. M. Christensen, \emph{Theory of Viscoelasticity,}  Acad. Press, New York,  1982.

 

\ZBI{Day} Day, W.A., On the monotonicity of the relaxation function of viscoelastic materials, \emph{Proc. Cambridge Philos. Soc. } 67,     503–508,  1970
  \ZBI{Dinzart} Dinzart, F., Lipi\'nski P., Self-consistent approach of the constitutive law of a two-phase visco-elastic material described by fractional derivative models. \emph{ Arch. Mech.}, {\bf 62}  135-156, 2010. 
\ZBI{FabrizioOWENS}Amendola, G., Fabrizio, M.,  Golden, J.M.,  \emph{Thermodynamics of materials with memory. Theory and applications.} Springer, New York, 2012. 
 
 

  \ZBI{Gerlac} Gerlac, S., Matzenmiller, A., Comparison of numerical methods for identification of viscoelastic line spectra from static tests. \emph{  Int. J. Numerical methods Eng.}, {\bf 63}  428-454, 2005.
  
  \bibitem{Gentili} G. Giorgi, G. Gentili, {  thermodynamic properties and stability for the heat flux equation with linear memory.\/} \emph{ Quarterly Appl. Math.}, ~{\bf 51} 343-362, 1993.
  
\ZBI{GolubKozbarRegulina} Golub V.P., Kobzar',  Yu. M., Ragulina V.S., Determining the parameters of the hereditary kernels in nonlinear viscoelastic materials in tension. \emph{ Int. Appl. Mechanics}, {\bf 49} 102-113,  2013. 

\ZBI{DEkEELIUHinestroza}   De Kee, D.   Liu Q. \&   Hinestroza, J. (2005). Viscoelastic (non-fickian) diffusion. {\em The Canada J. of chemical engineering,\/} {  83,}   913--929.


\ZBI{GRASSEKABANlore}  
Grasselli, M., Kabanikhin, S. I., Lorenzi, A., An inverse hyperbolic integrodifferential problem arising in geophysics. II. Nonlinear Anal. 15 (1990), no. 3, 283–298.

\ZBI{GRASSEKABAN}
Grasselli, M., Kabanikhin, S. I., Lorenzi, A.,  An inverse problem for an integro-differential equation.   \emph{Siberian Math. J.} {\bf  33} (1992), no. 3, 415–426 (1993) 

  \ZBI{Guidetti} Guidetti, D.,  Reconstruction of  variation convolution kernel in an abstract wave equation. \emph{Forum Math.} {\bf 6} 1129-1160, 2010. 

\ZBI{Janno} Janno, J., Identification of weakly singular relaxation kernels in three-dymensional viscoelasticity,  \emph{J. Math. Anal. Appl.,} 262, 133-159, 2001.

\ZBI{Kolsky} H. Kolsky, \emph{Stress waves in solids,}   Dover Publ., New York, 1963.
\ZBI{Koosis} Koosis, P., \emph{Introduction to $H_p$ spaces,} Cambridg U.P.,  Cambridge 1998.



\bibitem{Lerenzilibro} Lorenzi, A., \emph{An introduction to identification problems via functional analysis.} VSP, Utrecht,  2001.

\ZBI{BykovMILLISECONDS} Bykov, D.L., Kazakov, A.V., Konovalov, D.N., Me'lnikov, V.P., Osavchuk, A.N., Peleshko, V.A.,     Identification of the model of nonlinear viscoelasticity of filled polymer materials in millisecond time range. \emph{Mechanics of Solids}, {\bf 47}  641-645, 2012.     

\ZBI{PandLIBRO} Pandolfi, L.,
\emph{Distributed systems with persistent memory. 
Control and moment problems.} Springer Briefs in Electrical and Computer Engineering. Control, Automation and Robotics. Springer, Cham, 2014.

\ZBI{PandIDENT} Pandolfi, L., \emph{A linear algorithm for the identification of a relaxation kernel using two boundary measures,} Inverse Problems 31 (2015), no. 10, 105003, 12 pp.

\ZBI{PandMULTI} Pandolfi, L., \emph{Identification of the relaxation kernel in diffusion processes and viscoelasticity with memory via deconvolution,}  Math. Methods Appl. Sci., DOI: 10.1002/mma.4180
 
\ZBI{Pipkin} Pipkin, A.C., \emph{Lectures on viscoelasticity theory,} Springer, Berlin 1986.
\ZBI{Pruess} Pr\"uss, J., \emph{Evolutionary integral equations and applications,} Birkh\"auser Verlag, Basel, 1993. 
\end{thebibliography}
\enddocument